\documentclass[12pt]{amsart}
\usepackage{mathrsfs}
\usepackage{amsfonts}
\usepackage{amsmath}
\usepackage{amssymb}
\usepackage{graphicx}
\usepackage{color}
\usepackage{xcolor}
\usepackage{amsthm,amscd, amsrefs}
\usepackage{appendix,mathtools}
\usepackage{calc}
\usepackage{hyperref}
\usepackage{setspace}

\newtheorem{thm}{Theorem}[section]
\newtheorem{cor}[thm]{Corollary}
\newtheorem{lem}[thm]{Lemma}
\newtheorem{prop}[thm]{Proposition}

\newcommand{\R}{\mathbb R}

\newcommand{\n}{\nabla}

\newcommand{\vn}{\vec{n}}

\newcommand{\pd}{\partial}

\newtheorem{theorem}{Theorem}
\newtheorem{lemma}[theorem]{Lemma}

\newtheorem{definition}[theorem]{Definition}
\newtheorem{assumption}{Assumption}
\newtheorem{remark}[theorem]{Remark}

\newtheorem{conjecture}[theorem]{Conjecture}

\title[]{Upper bound of the counting function of Steklov eigenvalues}

\author{Fei He}
\email{hefei@xmu.edu.cn}
\address{School of Mathematical Science, Xiamen University, 422 S. Siming Rd. Xiamen, Fujian, P.R.China, 361005.}

\author{Lihan Wang}
\email{lihan.wang@csulb.edu}
\address{Department of Mathematics and Statistics, California State University, Long Beach, 1250 Bellflower Blvd, Long Beach, CA 90840.}

\thanks{
The first named author was partially supported by NSFC grant No.12141101.\\
	The second named author was partially supported by NSF grant DMS-2316620.}

\begin{document}

\begin{abstract}
We study the counting function of Steklov eigenvalues on compact manifolds with boundary and obtain its upper bound involving the leading term of Weyl's law. Our estimate can be viewed as a weakened version of P\'{o}lya's Conjecture in the Steklov case on general manifolds. As a byproduct, we also obtain a description about the decay behavior of Steklov eigenfunctions near the boundary.
\end{abstract}

\maketitle

\section{Introduction}

Given a smooth compact Riemannian manifold $M^n$ with smooth boundary $\partial M^n$, the Steklov  eigenvalue problem is defined as follows: 
\begin{align*}
\left\{ 
\begin{aligned}
\Delta u&=0,\, && M^n;\\ 
\frac{\partial u}{\partial \vec{n}}&=\sigma u,\, &&\partial M^n,
\end{aligned}%
\right.
\end{align*} where $\Delta$ is the Laplace-Beltrami operator on $M^n$ and $\vec{n}$ is the outward normal vector. Its spectrum is discrete and 
consists of isolated eigenvalues of finite multiplicity:
\begin{align*}
0=\sigma_0<\sigma_1\leq \sigma_2\leq \cdots\nearrow +\infty.
\end{align*}

This problem were first discussed by Steklov (\cite{S}) in 1902 motivated by physics: these eigenfunctions represent the steady state temperature on bounded domains such that the flux on the boundary is proportional to the temperature. The Steklov eigenvalue problem appears in many physical fields such as fluid mechanics, electromagnetism and elasticity. It has lots of applications in physics and technology like seismology and tomography. Mathematically, there has been intense interest and study in this problem. One recent breakthrough was made by Fraser and Schoen in the extremal Steklov eigenvalue problems on surfaces. In their seminal work \cite{FS16}, Fraser and Schoen revealed a deep connection between the extremal Steklov eigenvalue problems and the free boundary minimal surface theory in the unit Euclidean ball $\mathbb{B}^n$. See \cite{GP} for a review and \cite{CGGS} for  the recent development about Steklov eigenvalue problems.

We are interested in the counting function for Steklov eigenvalues :
\begin{align}\label{CF}
N(\sigma)&:=\# \{k\in \mathbb{N}: \sigma_k <\sigma\}.
\end{align} As shown by L. Sandgren in \cite{S55},  $N(\sigma)$ satisfies the Weyl's law when the boundary is $C^2$:
\begin{align}\label{Weyl}
N(\sigma)&=\frac{\omega_{n-1}}{(2\pi)^{n-1}}\rm{Vol}(\partial M)\sigma^{n-1}+o(\sigma^{n-1})
\end{align}with $\omega_{n-1}$ as the volume of the unit Euclidean ball $\mathbb{B}^{n-1}$. This asymptotic formula has been extended to domains with piecewise $C^1$ boundary in \cite{A06}, Lipschitz boundary in the two-dimensional case in \cite{KLP23} and higher dimensions in \cite{R23}. We refer the readers to \cite{KLP23} and \cite{R23} for more discussions and references about the asymptotic behavior of $N(\sigma)$. 

Our purpose in this paper is to estimate the counting function $N(\sigma)$ from above by the leading term of the Weyl's law \ref{Weyl}. This type of question has always been one of central interests in spectral geometry. In the case of the Laplace operator $\Delta$, such question is the well-known {\it P\'olya's Conjecture} (\cite{P54}, 1954) in spectral geometry:
\begin{conjecture}[P\'olya's Conjecture]
The eigenvalue counting functions of the Dirichlet Laplacian on a bounded Euclidean domain can be estimated from above by the leading term of Weyl’s Law.
 \end{conjecture}
Note the version of P\'olya's Conjecture for Neumann Laplacian states that its eigenvalue counting functions can be estimated from below by the leading term of the Weyl’s law. In \cite{P61}, Pólya  proved this conjecture for plane-covering domains. For general domains, Li and Yau proved a weakened version in \cite{LY83}. The P\'olya Conjecture for both of Dirichlet and Neumann Laplacian is still open in general. The recent breakthrough shows that it is true for Euclidean balls in \cite{FLPS23}. Please see \cite{Lin17}, \cite{FreS23} and \cite{FLPS23} for more results and recent development about the P\'olya's Conjecture.

In this paper, we proved the following weakened version of P\'olya's Conjecture for Steklov eigenvalues on compact manifolds with boundary:

\begin{thm}[Main Theorem]\label{intro - thm: main theorem}
Let $(M^n, g)$ be a smooth compact manifold with smooth boundary $\pd M^n$. Let $\rho$ be the distance function to the boundary and $\Sigma_\rho$ be its level sets. Let $\mathrm{I\!I}$ denote the second fundamental form and $H$ the mean curvature. Let $\rm{nir}(\partial M^n)$ denote the the normal injective radius of $\partial M^n$ and define
\begin{align*}
C_H=\sup_{\rho <\rm{nir}(\partial M)}\sup_{\Sigma_{\rho}} |H|,\quad C_{\mathrm{I\!I}}=\sup_{\rho <\rm{nir}(\partial M)}\sup_{\Sigma_{\rho}}(|\mathrm{I\!I}|).
\end{align*}
Assume that $Ric(x) \geq - (n-1)K$ for some constant $K \geq 0$ when $\rho(x) < \rm{nir}(\partial M)$. Then for any $\sigma>0$, there is \begin{align}\label{eqn:main thm1}
N(\sigma)&\leq C(n)e^{\left(C(n) \sqrt{K}+C_H\right)\rm{nir}(\partial M^n) + 4e^{6(C_H+C_{\mathrm{I\!I}} ) \min\{ {\rm nir}(\pd M) , \sigma^{-1}\} }} v_0^{-1} {\rm Vol}(\partial M) \left(\frac{1 }{\rm{nir}(\partial M^n)}+\sigma\right)^{n-1}.
\end{align} In particular, when $\sigma \geq \frac{c}{\rm{nir}(\partial M^n)}$ for any positive constant $c$,  there is 
\begin{align}\label{eqn:main thm2}
N(\sigma)&\leq C(n)e^{\left(C(n) \sqrt{K}+C_H\right)\rm{nir}(\partial M^n)+ 4e^{6c(C_H+C_{\mathrm{I\!I}} )\sigma^{-1}}} v_0^{-1} {\rm Vol}(\partial M) \sigma^{n-1}.
\end{align} 
The constant $v_0$ is the volume non-collapsing constant defined as \[v_0 = \inf_{\rho(x) < {\rm nir}(\pd M), 0< s< \min \{ {\rm nir}(\pd M)/2, (2\sigma)^{-1} \} } \frac{{\rm Vol}(B(x, s))}{s^n}.\] It can be replaced by a constant depending only on $n$ when $\sigma$ is sufficiently large. 
\end{thm}
The normal injective radius of the boundary is defined as the maximum of the set of radius $r>0$ such that the exponential map $\exp(x,t)=\exp_{x}(t \nu): \partial M \times [0, r)\rightarrow M$ is a diffeomorphism onto its image with $\nu$ as the inward unit normal of $\partial M$. In this theorem, the upper bounds are sharp in the sense that they have the same orders as in the Weyl's law \eqref{Weyl} in terms of $\rm{Vol}(\partial M)$ and $\sigma$. With the assumption of general compact manifolds instead of Euclidean domains, it is reasonable to expect a different coefficient from $\frac{\omega_{n-1}}{(2\pi)^{n-1}}$ in \eqref{Weyl} which only depends on the dimension. The coefficients in our estimates \ref{eqn:main thm1} and \ref{eqn:main thm2} depend on the geometry near the boundary including the second fundamental form, the mean curvature, the normal boundary injective radius and the lower bound of the Ricci curvature. 

The upper bound of $N(\sigma)$ we obtained has two immediate consequences about the Steklov eigenvalue $\sigma_k$ and its multiplicity.  Since $N(\sigma)$ is the number of eigenvalues $\sigma_k< \sigma$, its upper bound implies a lower bound of eigenvalues $\sigma_k$ as stated in Corollary \ref{cor: lower bound for eigenvalues}. In addition, $N(\sigma)$ is also the sum of the multiplicities of distinct eigenvalues less than $\sigma$. Therefore its upper bound implies the upper bound of the multiplicity of eigenvalues $\sigma_k$.  About the multiplicity of Steklov eigenvalues, the second author in \cite{W22} proved that non-zero Steklov eigenvalues are simple for generic Riemannian metrics.

To prove Theorem \ref{intro - thm: main theorem}, we follow the spirit of P.Li's work in estimating the dimension of the space of harmonic functions with polynomial growth in \cite{Li97}. Two key ingredients in P.Li's approach are the mean value inequality and the weak volume comparison condition. One change we make is to use the volume non-collapsing property of small balls, which is more suitable on compact manifold, instead of the weak volume comparison condition. A major challenge in implementing P.Li's idea in our case is to control the growth of Steklov eigenfunctions. To overcome this challenge, we adopt the method of frequency developed by the first named author and J. Ou in \cite{HeOu24} to estimate the growth of the $L^2$ norm of Steklov eigenfunctions. The frequency for harmonic functions is introduced by Almgren \cite{Alm77} and has been generalized and applied in various literature. We use the frequency to measure the growth of a weighted average of harmonic functions on level sets. The growth estimate implies a control of the growth of the $L^2$ norm of Steklov eigenfunctions on small balls (see Lemma \ref{L2 estimate on ball}). Detailed definitions and statements are given in Section \ref{section: frequency} and Section \ref{section: counting function in general case}. One observation which plays an important part in our proof is that the frequency of a Steklov eigenfunction can be naturally estimated from above using the corresponding eigenvalue.

We first carry out our approach in a more general situation in Section \ref{section: frequency} and Section \ref{section: counting function in general case}. Namely, we look at manifolds with a good potential function and estimate the counting function for Steklov eigenvalues on the sub-level sets of this potential function. This general setup actually covers many situations. As a first showcase we study the problem on star-shaped Euclidean domains in Section \ref{section: Euclidean domains}. Theorem \ref{intro - thm: main theorem} is proved in Section \ref{section: manifolds with smooth boundary}. In this case, we use the distance function to the the boundary to define a potential function. Then we refine the general estimates in Section \ref{section: frequency} and Section \ref{section: counting function in general case} to obtain the more explicit estimates in Theorem \ref{intro - thm: main theorem}.

As a byproduct we also get a description about the decay behavior of Steklov eigenfunctions near the boundary in $L^2$-integral form in Section \ref{section: manifolds with smooth boundary}.
\begin{thm}\label{intro - thm: behavior of eigenfunction near boundary}
Let $(M^n, g)$ be a smooth compact manifold with smooth boundary $\pd M^n$. Let $\rho$ be the distance function to the boundary and $\Sigma_\rho$ be its level sets. Let $H$ denote the mean curvature. Let $u$ be a Steklov eigenfunction with respect to eigenvalue $\sigma>0$. Then there is
\begin{align*}
 e^{(-2 \sigma -\inf H + o(1))\rho} \leq \frac{\int_{\Sigma_\rho}  u^2 }{ \int_{\partial M^n} u^2 } \leq   e^{(-2 \sigma - \sup H + o(1))\rho},
\end{align*} for sufficiently small $\rho$.
\end{thm}

This result is motivated by the commonly held hypothesis that Steklov eigenfunctions exhibit boundary-localized oscillations and rapid decay in the interior for large eigenvalues. In 2001, Hislop and Lutzer \cite{HL01} proved that the Steklov eigenfunctions decay super-polynomially into the interior of Euclidean domains when the boundary is smooth. And they conjectured that the decay should be exponential in the case of an analytic boundary. Galkowski and Toth \cite{GT19} proved this Hislop-Lutzer conjecture using point-wise upper bounds of Steklov eigenfunctions near the boundary. Recently, they also \cite{GT23} proved an $L^2$-lower bounds of Steklov eigenfunctions near the boundary for analytic boundary. Note that the $L^2-$upper and lower bounds we obtain in Theorem \ref{intro - thm: behavior of eigenfunction near boundary} exhibit exponential decay in eigenvalues and the distance-to-boundary function, and these estimates hold for manifolds with a smooth boundary.

 We would like to point out that the smoothness in Theorem \ref{intro - thm: main theorem} and Theorem \ref{intro - thm: behavior of eigenfunction near boundary} can be relaxed to $C^2$. In fact our proofs only need that the distance function to the boundary is $C^2$.

Finally in Section \ref{section: c1alpha domains} we study the case of $C^{1, \alpha}$ domains in Riemannian manifolds. By constructing a distance-like function to the boundary, we generalize Theorem \ref{intro - thm: main theorem} to $C^{1, \alpha}$ domains but with less explicit estimate. 
\begin{thm}\label{intro-thm: estimate of counting function - C1alpha domain}
Let $\Omega$ be a bounded $C^{1, \alpha}$ domain in a smooth Riemannian manifold $(M^n, g)$. There are constants $i_0$ and $C$ depending on the geometry of a neighborhood of $\Omega$ such that
\begin{align*}
&{N}(\sigma) \leq C Area(\pd \Omega)  \left( \frac{1}{ i_0} + \sigma \right)^{n-1}.
\end{align*}
\end{thm} We also obtain a description about the behavior of Steklov eigenfunctions near the boundary in this case. See Theorem \ref{thm:C1alpha} for details.

\textbf{Acknowledgement:} The second named author would like to thank
the Tianyuan Mathematical Center in Southeast China for its hospitality during her visit when this work was partially finished. 

\section{Frequency of harmonic functions on manifolds with potential functions}\label{section: frequency}
Let $(M^n, g, f)$ be a complete Riemannian manifold with a smooth positive function $f$. Here $f$ is called the potential function. \begin{definition}
Define a symmetric $2$-tensor $T_{ij}$ \begin{equation}\label{eqn: the tensor T}
 T_{ij}=\frac{1}{2} g_{ij} - \n_i\n_j f,
\end{equation}
and a function $S$ 

\begin{equation}\label{eqn: the function S}
 S=f -  |\n f|^2 .
\end{equation}
\end{definition}

\begin{remark}
In fact, the tensor $T$ and the function $S$ appear naturally on gradient shrinking Ricci solitons as Ricci curvature tensor and scalar curvature respectively. Let $(M^n, g, f)$ be a gradient shrinking Ricci soliton. Then $f$ satisfies
\[\n_i \n_j f + R_{ij} = \frac{1}{2}g_{ij}.\]
This implies that $T$ is just the Ricci curvature tensor in this case.

Let $Scal$ denote the scalar curvature of a gradient shrinking Ricci soliton. It is known that 
$$Scal + |\n f|^2 - f = constant.$$
By modifying $f$ by a constant, this formula implies that $S$ is just the scalar curvature in this case. See Chapter 1 in \cite{Cho07} for more details. 
\end{remark}

We introduce the following function related to $f$:
\begin{equation*}
b = 2 \sqrt{f}.
\end{equation*}
Since $f$ is positive, by direct calculations, it follows that
\[
\n b = \frac{\n f}{\sqrt{f}}, \quad \n\n b = \frac{\sqrt{f} \n\n f - \frac{1}{2\sqrt{f}} \n f \otimes \n f }{f}.
\]
By (\ref{eqn: the tensor T}) and (\ref{eqn: the function S}), the function $b$ satisfies 
\begin{equation}\label{eqn: properties of b}
\begin{cases}
|\n b|^2 = 1 - \frac{4S}{b^2},\\
b\Delta b =  n - |\n b|^2 - 2tr T,
\end{cases}
\end{equation}
where $trT = g^{ij} T_{ij}$ is the trace of $T$.

We are going to define the frequency of harmonic functions on level sets of the function $b$. And the following assumptions are in need:
\begin{assumption}\label{basic assumptions}
\begin{enumerate}
\item Assume there are constants $1 \leq R_0<  R$, such that each $r \in [R_0, R]$ is a regular value of $b$. 
\item Assume that the sublevel set $\{b \leq r\}$ is compact for each $r \in[R_0, R]$. 
\end{enumerate}
\end{assumption}

\begin{definition}\label{definition of freq}
Suppose that Assumption 1 holds on $(M^n, g, f)$. Consider a harmonic function $u$ on $\{b<R\}$. For each $r \in [R_0, R]$, define \[
I(r)= r^{1-n} \int_{b=r} u^2 |\n b|,
\] and \[
D(r) = r^{2-n} \int_{b= r} u \langle \n u, \vn \rangle\] with $\vn$ as the unit outward normal vector on $\{b = r\}$. Then $I(r)>0$ and the frequency of $u$ is defined as 
\[
U(r) = \frac{D(r)}{I(r)}.
\]
\end{definition}

Here $I(r)$ is the weighted average of $u$ on the level set of $b$ and 
the usage of $|\n b|$ as a weight follows \cite{CM21}. Since $u$ is a non-trivial harmonic function, the positivity of $I(r)$ at regular values $r$ of $b$ follows from the maximum principle. We also notice that 
\begin{equation}\label{def of D}
D(r) = r^{2-n} \int_{b< r} |\n u|^2
\end{equation} by the divergence theorem.

We want to estimate the growth of $U(r)$ and $I(r)$. At first we derive formulas of  $I^{\prime}(r)$ and $D^{\prime}(r)$.
\begin{lem}\label{lem: derivative of I} At a regular value $r$ of $b$, we have
\begin{align*}
I'(r) &= \frac{2D(r)}{r} + r^{-n} \int_{b=r} \left( \frac{4n S}{r^2} -2 tr T \right) \left(1-\frac{4S}{r^2}\right)^{-1} u^2|\n b|.
\end{align*} and 
\begin{align*}
( \ln I(r) )' &= \frac{2U(r)}{r} + r^{-n}I^{-1}(r) \int_{b=r} \left( \frac{4n S}{r^2} -2 tr T \right) \left(1-\frac{4S}{r^2}\right)^{-1} u^2|\n b|.
\end{align*}

\end{lem}
\begin{proof}
Let $\vn$ as the unit outward normal vector on $\{b = r\}$, that is, $\vn=\frac{\nabla b}{|\nabla b|}$. By the divergence theorem, it follows that
\begin{align*}
I(r) & = r^{1-n} \int_{b=r} u^2 |\n b|
=r^{1-n} \int_{b=r} u^2 \langle \n b, \vn \rangle \\
 &= \int_{b<r} \langle \n u^2 , \n b \rangle  b^{1-n} - (n-1) u^2 b^{-n} |\n b|^2 + u^2 b^{1-n} \Delta b. \\
\end{align*}
Plug (\ref{eqn: properties of b}) into this equation and it follows that
\begin{align*}
I(r) = & \int_{b< r} \langle \n u^2 , \n b \rangle  b^{1-n} - (n-1) u^2 b^{-n} |\n b|^2 + u^2 b^{-n} (n - |\n b|^2 - 2trT) \\
= & \int_{b< r} \langle \n u^2 , \n b \rangle  b^{1-n} + u^2 b^{-n} \left( \frac{4n S}{b^2} - 2 tr T \right).\\
\end{align*}
Take the derivative on the last equality. By the co-area formula, it follows that 
\begin{align*}
I^{\prime}(r)=&r^{1-n} \int_{b=r} \langle \n u^2, \vn \rangle + r^{-n}  \int_{b=r} \left( \frac{4n S}{r^2} - 2 tr T \right) \frac{ u^2}{|\n b|}\\
=&\frac{2D(r)}{r} + r^{-n}  \int_{b=r} \left( \frac{4n S}{r^2} - 2 tr T \right) \frac{ u^2}{|\n b|}.
\end{align*}
By \eqref{eqn: properties of b}, there is 
\begin{align*}
\frac{1}{|\n b| }&=(1-\frac{4S}{b^2})^{-1}|\nabla b|.
\end{align*} Applying this equality into the second term of above formula of $I^{\prime}(r)$ yields the conclusion. 
\end{proof}

\begin{lem}\label{lem: derivative of D} At a regular value $r$ of $b$, we have
\begin{equation*}
\begin{split}
D'(r) = & 2 r^{2-n }\int_{b= r} \left| \frac{\pd u}{\pd \vn }\right|^2 |\n b|^{-1} +  4 r^{-n} \int_{b=r}S\left(|\n u|^2 - 2 \left| \frac{\pd u}{\pd \vn }\right|^2 \right) |\n b|^{-1}\\
&+2r^{1-n} \int_{b< r}\left( -tr T |\n u|^2 + 2 T(\n u, \n u)\right).
\end{split}
\end{equation*}
\end{lem}

\begin{proof}
Take the derivative of D with respect to $r$ using \refeq{def of D}. It follows that 
\begin{equation}\label{deri of D}
D'(r) = r^{2-n}\left( \int_{b=r} \frac{|\n  u|^2}{ |\n b| }  - \frac{n-2}{r} \int_{b< r} |\n u|^2  \right).
\end{equation}
 To calculate the right hand side above, we need the following formula:
\begin{equation}\label{eqn: first variation formula}
\int_{b< r} |\n u|^2 div( X) - 2 \langle \n u \otimes \n u , \n X \rangle = \int_{b=r} |\n u|^2 \langle X, \vn \rangle - 2 \frac{\partial u}{\pd \vn} \langle \n u, X\rangle,
\end{equation}
for any smooth vector field $X$. This can be derived from the variation formula of the Dirichlet energy at $u$ w.r.t. deformations generated by $X$.  
It can also be proved directly using integration by parts. 

Take $X = \n f$ in (\ref{eqn: first variation formula}), we get
\begin{equation}\label{eqn: first variation formula with X = nabla f}
\begin{split}
\int_{b< r} |\n u| ^2 \Delta f - 2 f_{ij}u_i u_j = & \int_{b=r} |\n u|^2 |\n f| - 2 \langle \n u, \vn \rangle \langle \n u, \n f\rangle  \\
= & \frac{r}{2} \int_{b=r} |\n u|^2 |\n b| - 2 \left| \frac{\pd u}{\pd \vn} \right|^2 |\n b|.
\end{split}
\end{equation}
Since $|\n b|^2 = 1 - \frac{4S}{b^2}$ by \eqref{eqn: properties of b}, the RHS of \eqref{eqn: first variation formula with X = nabla f} can be written as
\begin{equation}\label{RHS}
RHS = \frac{r}{2} \int_{b= r} \frac{|\n u|^2 }{|\n b|} \left( 1 - \frac{4S}{b^2}\right) - 2\frac{ \left| \frac{\pd u}{\pd \vn} \right|^2}{|\n b|} \left( 1 - \frac{4S}{b^2}\right),
\end{equation}
By \eqref{eqn: the tensor T}, we have

\[
\Delta f = \frac{n}{2} - tr T, \quad f_{ij} = \frac{1}{2} g_{ij} - T_{ij}.
\]
Plug these into the LHS of (\ref{eqn: first variation formula with X = nabla f}) and then it follows that 

\begin{equation}\label{LHS}
LHS = \int_{b< r} \frac{n-2}{2} |\n u|^2 - tr T |\n u|^2 + 2 T (\n u, \n u).
\end{equation}

Then the equality of \eqref{LHS} and \eqref{RHS} implies that 
\begin{align*}
\frac{r}{2} \int_{b= r} \frac{|\n u|^2 }{|\n b|}-\int_{b< r} \frac{n-2}{2} |\n u|^2&=\int_{b= r} \left(\frac{2}{r} S|\n u|^2  + r \left| \frac{\pd u}{\pd \vn} \right|^2  -\frac{4}{r}S \left| \frac{\pd u}{\pd \vn} \right|^2\right)|\n b|^{-1}\\
&\int_{b< r} \left(-tr T |\n u|^2 + 2 T (\n u, \n u)\right).
\end{align*}
Plug this equality into the right hand of \eqref{deri of D} after multiplying $2r^{1-n}$ and then desired formula follows.
\end{proof}

In order to estimate $U^{\prime}(r)$ from below, we need a lower bound of $D^{\prime}(r)$.

\begin{lem}\label{lem: iterating the first variation formula}
Assume that $\frac{ 4 }{r^2}\sup_{b< r}|S|<1$ for all $R_0 \leq r\leq R$. Then there is a constant $C$ depending on $\sup_{R_0 \leq r\leq R}\left(\frac{ 4 }{r^2}\sup_{b<r}|S|\right)<1$ and the dimension $n$, such that

\begin{align*}
D'(r) \geq 2 r^{-1 }D^2(r)I^{-1}(r) - \frac{C\sup_{b< r} (|\n S||\n b| + |S\n\n b| + |T|) }{r} D(r)
\end{align*}
for $1<R_0\leq r \leq R$. 

\end{lem}
\begin{proof}
By Lemma \ref{lem: derivative of D}, there is 
\begin{equation}\label{deri of D2}
\begin{split}
D'(r) = & 2 r^{2-n }\int_{b= r} \left| \frac{\pd u}{\pd \vn }\right|^2 |\n b|^{-1} +  4 r^{-n} \int_{b=r}S\left(|\n u|^2 - 2 \left| \frac{\pd u}{\pd \vn }\right|^2 \right) |\n b|^{-1}\\
&+2r^{1-n} \int_{b< r}\left( -tr T |\n u|^2 + 2 T(\n u, \n u)\right).
\end{split}
\end{equation}
We first consider the first term on the right hand side of \eqref{deri of D2}. By the definition of $D(r)$ and $I(r)$, we have 
\begin{align*}
D^2(r)&=r^{4-2n}\left|\int_{b=r} u\frac{\pd u}{\pd \vn }\right|^2\\
&\leq r^{4-2n}\int_{b=r}u^2|\nabla b|\int_{b=r}\left|\frac{\pd u}{\pd \vn }\right|^2|\nabla b|^{-1}\\
&=r^{3-n}I(r)\int_{b=r}\left|\frac{\pd u}{\pd \vn }\right|^2|\nabla b|^{-1}.
\end{align*} Here Holder's inequality is applied in the second line. Then we have
\begin{equation}\label{deri of D3}
\begin{split}
D'(r) \geq & 2 r^{-1 }D^2(r)I^{-1}(r)+  4 r^{-n} \int_{b=r}S\left(|\n u|^2 - 2 \left| \frac{\pd u}{\pd \vn }\right|^2 \right) |\n b|^{-1}\\
&+2r^{1-n} \int_{b< r}\left( -tr T |\n u|^2 + 2 T(\n u, \n u)\right).
\end{split}
\end{equation}

Next we estimate the second term from below. For any positive integer $m$ and $R_0\leq r \leq R$, define
\[
K_m(r) = \int_{b=r} S^m \left( |\n u|^2 - 2 \left| \frac{\pd u}{\pd \vn }\right|^2 \right) |\n b|^{-1}.
\] Then $4r^{-n}K_1(r)$ is exactly the second term in \eqref{deri of D2}.

By \eqref{eqn: properties of b}, we have the identity $|\nabla b|^{-1} = |\n b| + \frac{4S}{b^2|\nabla b|}$. Plugging this into the integrand above yields that
\[
\begin{split}
K_m(r) = \int _{b=r}  S^m\left( |\n u|^2 - 2 \left| \frac{\pd u}{\pd \vn }\right|^2 \right) |\n b| + \frac{4}{r^2}K_{m+1}(r).
\end{split}
\]
By applying (\ref{eqn: first variation formula}) with $X = S^m \n b$, we have
\begin{equation*}\label{eqn: induction formula for K_m}
\begin{split}
& K_m(r) - \frac{4}{r^2} K_{m+1}(r) \\
= &  \int_{b< r} |\n u|^2 div\left( S^m \n b\right) - 2 \langle \n u \otimes \n u , \n (S^m \n b) \rangle \\
= & \int_{b< r}  S^{ m-1 } \left( |\n u|^2 (m \langle \n S, \n b\rangle  + S \Delta b )- 2  \langle \n u \otimes \n u , m \n S \otimes \n b + S \n\n b \rangle \right)\\
\geq & - m C_1 (\underset{b<r}{\sup}|S|)^{m-1} \int_{b< r} |\n u|^2.
\end{split}
\end{equation*}
Here we use the notation $C_1 = \sup_{b < r} (3|\n S| |\n b| + (n+ 2) |S| |\n\n b|)$, which is independent of $m$, for convenience.  

By iterating the above inequality, we get
\begin{align}\label{eqn: iterating K_m}
\begin{aligned}
K_1(r) &\geq \left( \frac{4}{r^2} \right)^m K_{m+1}(r)- \sum_{k=1}^m k \left( \frac{ 4 }{r^2}\underset{b<r}{\sup}|S| \right)^{k-1}  C_1 \int_{b< r} |\n u|^2\\
&\geq \left( \frac{4}{r^2} \right)^m K_{m+1}(r)- \sum_{k=1}^m k \left( \sup_{R_0 \leq r\leq R}\left(\frac{ 4 }{r^2}\sup_{b<r}|S|\right) \right)^{k-1}  C_1 \int_{b< r} |\n u|^2.
\end{aligned}
\end{align}
Since $\frac{ 4 }{r^2}\sup_{b<r}|S|< 1$ by the assumption, there is 
\begin{align*}
\left|\left( \frac{4}{r^2} \right)^m K_{m+1}(r)\right |&\leq \left( \frac{4}{r^2} \underset{r<R}{\sup} |S|\right)^m \int_{b=r}\left| |\n u|^2 - 2 \left| \frac{\pd u}{\pd \vn }\right|^2\right| \rightarrow 0 \, \text{as $m\rightarrow \infty$}.
 \end{align*}  And the series $\sum_{k=1}^{\infty} k \left(\sup_{R_0 \leq r\leq R}\left(\frac{ 4 }{r^2}\sup_{b<r}|S|\right)\right)^{k-1} $ is convergent to some constant $C$ depending on $\sup_{R_0 \leq r\leq R}\left(\frac{ 4 }{r^2}\sup_{b<r}|S|\right)<1$. Thus taking $m \to \infty$ in \eqref{eqn: iterating K_m} yields that 
\begin{align}\label{eqn:K_1}
K_1(r) &\geq - C C_1 \int_{b< r} |\n u|^2 \geq - C C_1r^{n-2}D(r)\geq - C C_1r^{n-1}D(r).
\end{align} The fact $r\geq R_0>1$ is used for the last inequality.

For the third integral term in \eqref{deri of D3}, it follows that
\begin{align}\label{eqn:3rd term}
\int_{b< r}\left( -tr T |\n u|^2 + 2 T(\n u, \n u)\right)&\geq - 3\sqrt{n} \sup_{\{b<r\}}|T| \int_{b<r}|\nabla u|^2= -3\sqrt{n} \sup_{\{b<r\}}|T|r^{n-2} D(r),
\end{align}
where we have used $|tr T| \leq \sqrt{n} |T|$.
Plug both of \eqref{eqn:K_1} and \eqref{eqn:3rd term} into \eqref{deri of D3}. Then the desired claim follows. 
\end{proof}

Next we will derive a lower bound of $U^{\prime}(r)$ which will control the growth of $U(r)$ and $I(r)$.
\begin{prop}\label{prop: differential inequality for U}
 Assume that $\sup_{R_0 \leq r \leq  R}\left(\frac{ 4 }{r^2}\sup_{b<r}|S|\right)<1$. Then there is a constant $C$ depending on $\sup_{R_0 \leq r\leq R}\left(\frac{ 4 }{r^2}\sup_{b<r}|S|\right)$ and the dimension $n$, such that
\begin{align*}
U'(r) &\geq - \frac{C\sup_{b\leq r} (|\n S||\n b| + |S\n\n b| + |T|)}{r} U(r)
\end{align*}
for $1<R_0\leq r\leq R$.

\end{prop}

\begin{proof}
By the definition of $U$, we have
\begin{align}\label{eqn:IU}
 I^2(r)U'(r) &= D'I - D I'.
\end{align} 

By Lemma \ref{lem: derivative of I}, there is
\begin{align}\label{eqn: upper of I^{prime}}
\begin{aligned}
I'(r)&\leq \frac{2D(r)}{r} + \underset{b=r}{\sup} \left[\left( \frac{4n S}{r^2} - 2 tr T \right) \left(1-\frac{4S}{r^2}\right)^{-1}\right] r^{-n}\int_{b=r}  u^2|\n b|\\
&=\frac{2D(r)}{r}  +\underset{b=r}{\sup} \left[\left( \frac{4n S}{r^2} -2 tr T \right) \left(1-\frac{4S}{r^2}\right)^{-1}\right] r^{-1}I(r).
\end{aligned}
\end{align}
Apply this inequality \eqref{eqn: upper of I^{prime}}
and Lemma \ref{lem: iterating the first variation formula} to \eqref{eqn:IU}. Then it follows
\begin{align*}
I^2(r)U'(r) & \geq -\frac{C\sup_{b< r} (|\n S||\n b| + |S\n\n b| + |T|) }{r} D(r)I(r)\\
 &-\sup_{b=r} \left[\left( \frac{4n S}{r^2} - 2 tr T \right) \left(1-\frac{4S}{r^2}\right)^{-1}\right] r^{-1}D(r)I(r).
\end{align*}
Since 
\begin{align*}
|\left( \frac{4n S}{r^2} - 2 tr T \right) \left(1-\frac{4S}{r^2}\right)^{-1}|&\leq \left(\frac{ 4 }{r^2}\sup_{b<r}|S|+ 2 \sqrt{n}\sup_{b\leq r} |T|\right) \left(1-\frac{ 4 }{r^2}\sup_{b<r}|S|\right)^{-1}\\
&\leq C(1+\sup_{b\leq r} |T|)
\end{align*} for some constant $C$ depending on $\sup_{R_0< r< R} (\frac{4}{r^2}\sup_{b<r}|S|)$ and the dimension $n$, there is 
\begin{align*}
I^2(r)U'(r) & \geq -\frac{C\sup_{b \leq r} (|\n S||\n b| + |S\n\n b| + |T|) }{r} D(r)I(r).
\end{align*} Therefore the conclusion follows since $U(r)=\frac{D(r)}{I(r)}$.
\end{proof}
Now we are ready to prove the following estimate on the growth of $U(r)$ and $I(r)$.
\begin{theorem}\label{thm: upper bd UI}
Assume that $\sup_{R_0 \leq r\leq R}\left(\frac{ 4 }{r^2}\sup_{b<r}|S|\right)< 1$. Assume there exists a constant $\kappa$ and $\beta \in [0, 1)$ such that 
 \begin{align*}\sup_{b< r} (|\n S||\n b| + |S\n\n b| + |T|) \leq \kappa (R - r)^{-\beta}.
 \end{align*}
Then there is a constant $C$ depending on $\sup_{R_0 \leq r\leq R}\left(\frac{ 4 }{r^2}\sup_{b<r}|S|\right)$ and $n$, such that
\begin{align}\label{eqn:upper bd of U}
U(r) & \leq \kappa C U(R)\exp\left(  \int_{R_0} ^ R \frac{ds}{s(R-s)^\beta} \right).
\end{align} And for $R_0\leq r_2<r_1\leq R$, there is

\begin{align}\label{eqn:upper bd of I}
\left(\frac{r_1}{r_2}\right)^{-C_2} I(r_2)\leq  I(r_1) &\leq \left(\frac{r_1}{r_2}\right)^{\left(C\kappa U(R)\exp\left(  \int_{R_0} ^ R \frac{ds}{s(R-s)^\beta}\right) + C_2\right)}I(r_2)
\end{align} 
with 
\begin{align*} 
C_2 &= \sup_{R_0 \leq b \leq R}  \left[ (\frac{4n S}{b^2} - 2 tr T)(1 - \frac{4S}{b^2})^{-1} \right].
\end{align*}

\end{theorem}
\begin{proof}
By Proposition \ref{prop: differential inequality for U}, we have
\[
U'(r) \geq - \frac{\kappa C}{r(R- r)^\beta} U(r)
\]
where $C$ depends on $\sup_{R_0 \leq r\leq R}(\frac{4}{r^2}\sup_{b<r}|S|)$ and $n$. Integrating the above inequality yields
\begin{align*}
U(r) &\leq \kappa CU(R) \exp\left(\int_r ^ R \frac{ds}{s(R-s)^\beta}  \right) \leq \kappa C U(R)\exp\left(  \int_{R_0} ^ R \frac{ds}{s(R-s)^\beta} \right).
\end{align*}

For convenience, let $C_3=\kappa CU(R) \exp\left(\int_{R_0} ^ R \frac{ds}{s(R-s)^\beta}  \right)$. Applying the above upper bound of $U(r)$ to the formula of $(\ln I(r) )'$ in Lemma \ref{lem: derivative of I}, we have
\begin{align*}
-\frac{C_2}{r} \leq ( \ln I(r) )'  &\leq \frac{2C_3+C_2}{r}
\end{align*}
with $C_2 = \sup_{R_0 \leq b \leq R}  \left[ (\frac{4n S}{b^2} - 2 tr T)(1 - \frac{4S}{b^2})^{-1} \right]$. Then integrating the above inequalities of $( \ln I(r) )'$ over $[r_2, r_1]\subset [R_0, R]$ implies  \eqref{eqn:upper bd of I}.

\end{proof}

\section{Counting functions on manifolds with potential functions}\label{section: counting function in general case}

With the set-up in Section \ref{section: frequency}, we will derive a general upper bound of the counting function $N(\sigma)$ for Steklov eigenvalues on $(M^n, g, f)$. Consider the Steklov eigenvalue problem on $\{b<R\}$:
\begin{align}\label{eqn:Steklov1}
\begin{cases}
\Delta u = 0 & \text{on } \{b<R\}; \\
\frac{\pd u}{\pd \vn} = \sigma u & \text{on } \{b=R\}.
\end{cases}
\end{align}
Through out this article, we only consider the non-trivial case, i.e., $\sigma \neq 0$ or $u$ is non-constant. Let $E_\sigma$ be the linear space spanned by all Steklov eigenfunctions with eigenvalue $\leq \sigma$.  Then the dimension of $E_\sigma$ is exactly the counting function $N(\sigma)$. 

We are going to apply Section 1 to functions in $E_\sigma$ and obtain the upper bound of the dimension of $E_{\sigma}$, i.e., $N(\sigma)$. First we have the following observation about the frequency of functions in $E_\sigma$.
\begin{lemma}
For any $u \in E_\sigma$, its frequency satisfies 
\begin{align}\label{eqn: frequency upper bound on boundary}
U(R) \leq \frac{\sigma R}{\inf_{b = R} |\n b|}.
\end{align}
\end{lemma}
\begin{proof}
Let $k$ denote the dimension of $V_{\sigma}$. Then there exists Steklov eigenfunctions $u_1, ...,u_k $ corresponding to eigenvalues $\sigma_1, ..., \sigma_k$ in \eqref{eqn:Steklov1} which form an orthonormal basis for $V_\sigma$ w.r.t. the $L^2$ inner product on $\{b=R\}$.

For any $u \in E_{\sigma}$, there exists constants $l_1, \cdots, l_k$ such that $u = \sum_{i = 1}^k l_i u_i$. By the definition \ref{definition of freq}, its frequency satisfies that:
\begin{align*}
U(R) = \frac{ R \int_{b = R} u \frac{\pd u}{\pd \vn}}{ \int_{b = R} u^2 |\n b|} = \frac{R\int_{b = R} \sum l_i^2 \sigma_i u_i^2}{ \int_{b = R} u^2 |\n b|}\leq \frac{\sigma R}{\inf_{b = R} |\n b|}.
\end{align*}
\end{proof}

\begin{thm}\label{thm: estimate of spectral counting function - general case}
Suppose Assumption \ref{basic assumptions}  and $\sup_{R_0 \leq r\leq R}\left(\frac{ 4 }{r^2}\sup_{b<r}|S|\right)<1$ hold.  Assume  there exists a constant $\kappa$ and $\beta \in [0, 1)$ such that 
 \begin{align*}\sup_{b< r} (|\n S||\n b| + |S\n\n b| + |T|) \leq \kappa (R - r)^{-\beta}.
 \end{align*}
 Let $R_0 =\alpha R$ for some $\alpha \in (0,1)$. Let $\gamma_0 = \inf_{R_0<b < R} |\n b|$ and $\gamma_1 = \sup_{R_0< b < R} |\n b|$. Then 
\begin{align}\label{eqn: bd of N}
{N}(\sigma) &\leq C_M C_7 A(R) v_0^{-1}\left[\left([\alpha^{-1/3} - 1]R\right)^{-1}+\sigma\right]^{n-1}.
\end{align} 
Here $C_M$ denotes the constant in the mean value inequality for harmonic functions on small balls and depends on the Ricci lower bound and $n, R, \gamma_1$. The constant $C_7$ depends on $n, \alpha, \beta,\gamma_1, \gamma_0, \kappa, \sup_{R_0 \leq r\leq R}\left(\frac{ 4 }{r^2}\sup_{b<r}|S|\right)$ and $C_2=\sup_{R_0 \leq b \leq R}  \left[ (\frac{4n S}{b^2} - 2 tr T)(1 - \frac{4S}{b^2})^{-1} \right]$. The constant $A(R)=\sup_{R_0 < s <  R}A(\{b = s\})$ is the supremum of the area of level sets of $b$. And $v_0$ is a volume non-collapsing constant defined by 
$$v_0 = \inf_{R_0 < b(x)< R, s < \min \{ R-R_0, \sigma^{-1} \}/(2\gamma_1) }  s^{-n} {\rm Vol}(B(x, s)).$$
In particular, when $\sigma$ is sufficiently large, the constant $v_0$ can be replaced by a constant depending only on $n$.
\end{thm}

To prove this theorem, we need the following key lemma about the growth of $L^2$ norm of harmonic functions.

\begin{lemma}[Key lemma]\label{L2 estimate on ball}
With the assumptions in Theorem \ref{thm: estimate of spectral counting function - general case}, there is a constant $C$ depending on $\sup_{R_0 \leq r\leq R}\left(\frac{ 4 }{r^2}\sup_{b<r}|S|\right)$ and $n$, suppose the constant $\lambda > 1$ is properly chosen such that $\lambda^{-3}R > R_0$, then for each point $x \in \{ \lambda^{-2} R < b < \lambda^{-1} R \}$ such that $B(x, \delta) \subset \{\lambda^{-3}R< b< R\}$ for some $\delta > 0$, there is
\begin{align}\label{eqn: L2 estimate on ball}
 \int_{B(x, \delta)}u^2 
 &\leq  \gamma^2_1\gamma_0^{-2}(\lambda^{C_2 - n} + 1 + \lambda^{n+C_2 + C_3}) \int_{\lambda^{-2}R < b < \lambda^{-1}R} u^2 
\end{align} where 
 $C_2 = \sup_{R_0 \leq b \leq R}  \left[ (\frac{4n S}{b^2} - 2 tr T)(1 - \frac{4S}{b^2})^{-1} \right]$, and $C_3=\kappa C\frac{\sigma R}{\gamma_0} \exp\left(  \int_{R_0} ^ R \frac{ds}{s(R-s)^\beta} \right)$.

\end{lemma}

\begin{proof}[Proof of Lemma \ref{L2 estimate on ball}]
Choose any $u \in V_{\sigma}$ and consider its weighted $L^2$ norm $I(r)$ for $r \in [R_0, R]$. Applying \eqref{eqn: frequency upper bound on boundary} to \eqref{eqn:upper bd of I}, we have 
\begin{align*}
\left(\frac{r_1}{r_2}\right)^{-C_2} I(r_2)\leq  I(r_1) &\leq \left(\frac{r_1}{r_2}\right)^{\left(\kappa CU(R)\exp\left(  \int_{R_0} ^ R \frac{ds}{s(R-s)^\beta}\right) + C_2\right)}I(r_2)
\end{align*} 
Then for any constant $\lambda > 1$ satisfying $\lambda r \leq R$, the above inequalities imply that
\begin{align}\label{eqn: bound I}
\lambda^{-C_2} I(r)\leq  I(\lambda r) &\leq  \lambda^{(C_3 + C_2)}I(r).
\end{align}

Define
\begin{align*}
J(r_1, r_2) &= \int^{r_2}_{r_1}r^{n-1}I(r) dr.
\end{align*}
With a change of variable $ r = \lambda s$, we have
\begin{align*}
J(\lambda r_1, \lambda r_2) &= \int^{\lambda r_2}_{\lambda r_1}r^{n-1}I(r) dr= \lambda^n \int_{r_1}^{r_2} s^{n-1} I(\lambda s) ds.
\end{align*}Applying \eqref{eqn: bound I}to $I(\lambda s)$ implies the following estimates on annulus regions:
\begin{align*}
\lambda^{n - C_2} J(r_1, r_2) \leq J(\lambda r_1, \lambda r_2)&  \leq \lambda^{n + C_2 +  C_3} J(r_1, r_2).
\end{align*}
Then these estimates implies that 
\begin{equation}\label{eqn: doubling estimate on annulus regions}
\begin{split}
J(\lambda^{-3} R, R) = & J(\lambda^{-3} R, \lambda^{-2}R) + J( \lambda^{-2} R, \lambda^{-1}R ) + J(\lambda^{-1} R, R) \\
\leq & (\lambda^{C_2 - n} + 1 + \lambda^{n+C_2 +  C_3}) J( \lambda^{-2} R, \lambda^{-1} R ).
\end{split}
\end{equation}
On the other hand, by the co-area formula, there is
\begin{align}\label{eqn:J}
J(r_1, r_2)& = \int^{r_2}_{r_1}\left(\int_{b=r}u^2|\nabla b|\right) dr=\int_{r_1 < b < r_2} u^2 |\n b|^2. 
\end{align} Together with \eqref{eqn: doubling estimate on annulus regions}, this implies that
\begin{align}\label{eqn: L2 estimate}
\begin{aligned}
\int_{\lambda^{-3} R< b < R} u^2&\leq \gamma^{-2}_0\int_{\lambda^{-3} R< b < R} u^2 |\n b|^2= \gamma^{-2}_0J(\lambda^{-3} R, R)\\
&\leq \gamma^{-2}_0(\lambda^{C_2 - n} + 1 + \lambda^{n+C_2 +  C_3}) J( \lambda^{-2} R, \lambda^{-1} R ).
\end{aligned}
\end{align}

For each point $x \in \{ \lambda^{-2} R < b < \lambda^{-1} R \}$ such that $B(x, \delta)$ is contained in $\{\lambda^{-3} R< b < R\}$. By \eqref{eqn: L2 estimate}, we have
\begin{align*}
 \int_{B(x, \delta)} u^2 &\leq  \int_{\lambda^{-3} R< b < R} u^2\leq \gamma^{-2}_0(\lambda^{C_2 - n} + 1 + \lambda^{n+C_2 +  C_3}) J( \lambda^{-2} R, \lambda^{-1} R ).
\end{align*}At the same time \eqref{eqn:J} implies that
\begin{align*}
J( \lambda^{-2} R, \lambda^{-1}R)&\leq \gamma_1^2 \int_{\lambda^{-2}R < b < \lambda^{-1}R} u^2.
\end{align*}
 Therefore the conclusion \eqref{eqn: L2 estimate on ball} follows.

\end{proof}

\begin{proof}[Proof of Theorem \ref{thm: estimate of spectral counting function - general case}
]
Choose
\begin{align}\label{eqn: choice of lambda and delta}
\lambda &=1+ \frac{1}{(\alpha^{-1/3} - 1)^{-1} +R\sigma}
\end{align} 
Note that this choices of $\lambda$ implies $\lambda\leq \alpha^{-1/3}$ to make sure that $ R_0=\alpha R < \lambda^{-3} R$.

 We introduce the following inner product  on $E_\sigma$:
\begin{align*}
L(u, v) &= \int_{\lambda^{-2} R < b < \lambda^{-1}R} uv.
\end{align*}
Note that $L(u, u) = 0$ implies that $u \equiv 0$ on $\{b<R\}$ by the unique continuation principle of harmonic functions.

 Let $\{u_1, u_2, ..., u_k\}$ be an orthonormal basis of $E_\sigma$ with respect to $L(\cdot, \cdot)$ and $F(x)=\sum^{k}_{i=1}u_i^2(x)$. Then $F$ is well-defined under an orthonormal change of basis and 
 \begin{align*}
 k=\rm{dim} E_{\sigma}= \int_{\lambda^{-2} R < b < \lambda^{-1}R} F(x)
 \end{align*}
 
For each point $p \in \{ \lambda^{-2} R < b < \lambda^{-1} R \}$, let $E_p=\{u\in E_{\sigma}| u(p)= 0\}$. By an argument due to P. Li (see Lemma 7.3 in \cite{Li12}), the co-dimension of the subspace $E_p$ in $E_{\sigma}$ is at most one. For completeness, we include the argument here. Suppose there are at least two linearly independent functions $w_1$ and $w_2$ in the complement of $E_p$. Then $w_1(p) \neq 0$ and $w_2(p)\neq 0$. On the other hand, their linear combination $w_1(p)w_2-w_2(p)w_1$ vanishes at the point $p$. Hence $w_1(p)w_2-w_2(p)w_1 \in E_p$ and we get a contradiction.

Therefore by a change of orthonormal basis, for the point $p$, we can have an orthonormal basis of $E_\sigma$, still denoted by $\{u_1, u_2, ..., u_k\}$, such that
\begin{align*}
F(p) = u_k^2(p).
\end{align*}
Apply the well-known Li-Schoen's mean value inequality (see Theorem 7.2 in \cite{Li12}) to $u_k^2$ on $B(p, \delta)$ with $\delta=\frac{R}{2\gamma_1} (1-\lambda^{-1})\lambda^{-2} $. It follows that 
\begin{align}\label{eqn: F1}
F(p)=u_k^2(p) \leq C_4\frac{1+\exp (C_5 \delta \sqrt{K})}{V(B(p, \delta))} \int_{B(p, \delta)} u_k^2
\end{align} for constants $C_4(n), C_5(n)>0$. Here $K$ is from $Ric \geq -(n-1)K$ on $\{R_0 < b<R\}$. Since $\delta<\frac{R}{\gamma_1}$, for convenience, let $C_M=C_4(1+\exp (C_5 R\gamma^{-1}_1 \sqrt{K}))$ be the constant of the mean value inequality which depends on $n, K, R, \gamma_1$.

With the choice of $\delta=\frac{R}{2\gamma_1} (1-\lambda^{-1})\lambda^{-2} $, we have that $B(p, 2\delta)$ is contained in $\{b<R\}$ for every $p \in \{ \lambda^{-2} R < b < \lambda^{-1} R\}$. We also have $\delta < \min\{ \frac{R-R_0}{2\gamma_1}, \frac{1}{2\gamma_1 \sigma}\}$ by a simple calculation, hence 
\[
V(B(p, \delta)) > v_0 \delta^n
\]
by the assumption. Morever, when $\sigma$ is large enough such that $\delta$ is less than the half of the injective radius of $\{R_0 \leq b\leq R\}$, a well-known result of Croke (see Proposition 14 in \cite{Croke80}) yields that the constant $v_0$ depends only on the dimension $n$. We call $v_0$ the volume non-collapsing constant for small balls. Then \eqref{eqn: F1} becomes
\begin{align}\label{eqn: F2}
F(p) \leq \frac{C_M}{v_0 \delta^n} \int_{B(p, \delta)} u_k^2.
\end{align} 

Applying Lemma \ref{L2 estimate on ball} to the right hand side of \eqref{eqn: F2} implies that   
\begin{align}\label{eqn: F3}
F(p)  
  &\leq \frac{C_M}{ v_0 \delta^n}  \gamma^2_1\gamma_0^{-2}(\lambda^{C_2 - n} + 1 + \lambda^{n+C_2 + C_3})
\end{align}
 for any $p \in \{ \lambda^{-2} R < b < \lambda^{-1} R\}$ since $\int_{\lambda^{-2}R < b < \lambda^{-1}R} u_k^2=1$.
Integrating \eqref{eqn: F3} on $\{\lambda^{-2} R < b < \lambda^{-1} R \}$, we get
\begin{equation}\label{eqn: rough estimate of k}
k \leq  V({ \lambda^{-2} R < b < \lambda^{-1} R }) \frac{C_M}{ v_0 \delta^n}\gamma^2_1\gamma_0^{-2}(\lambda^{C_2 - n} + 1 + \lambda^{n+C_2 + C_3}).
\end{equation}

At the same time, we have
\begin{align*}
V( \lambda^{-2} R < b < \lambda^{-1} R) &=\int_{ \lambda^{-2} R}^{ \lambda^{-1} R}\int_{b = s} \frac{1}{|\n b|} \\
&\leq \frac{R(\lambda -1)}{\lambda^2}\frac{ \sup_{\lambda^{-2} R < s < \lambda^{-1} R}A(\{b = s\})}{\gamma_0}\\
&\leq \frac{R(\lambda -1)}{\lambda^2}\frac{A(R)}{\gamma_0}
\end{align*} with $A(R)=\sup_{R_0 < s <  R}A(\{b = s\})$ as the supremum of the area of level sets of $b$. Plug this volume estimate into \eqref{eqn: rough estimate of k}. Then by $\delta=\frac{R}{\gamma_1} (1-\lambda^{-1})\lambda^{-2}$, it follows 
\begin{align}\label{eqn: 2rough estimate of k}
k &\leq \frac{C_M\gamma_1^{n+2}A(R)}{v_0\gamma_0^3}\lambda^{C_6(C_2, C_3, n)} \left((\lambda-1)R\right)^{1-n}
\end{align} with the constant $C_6$ depending on $C_2, C_3, n$.

Plug the definition \eqref{eqn: choice of lambda and delta} of $\lambda$  into \eqref{eqn: 2rough estimate of k}. Then we have 
\begin{align*}
k &\leq C_M C_7A(R)v_0^{-1} \left(\left((\alpha^{-1/3} - 1)R\right)^{-1}+\sigma\right)^{n-1}
\end{align*} with $C_7$ depending on $n, \gamma_1, \gamma_0, C_2, C_3, \alpha, \beta$.

Since $\int_{\alpha R} ^ R \frac{ds}{s(R-s)^\beta}  \leq \alpha^{-1}R^{-\beta}(1-\beta)^{-1}$, there is 
\begin{align*}
C_3 &\leq  \kappa C \frac{\sigma R}{\gamma_0} \exp\left( \frac{\alpha^{-1}R^{-\beta}}{1-\beta} \right)
\end{align*} where $C$ depending on $\sup_{R_0 \leq b\leq R}\left(\frac{4}{r^2}\sup_{b<r} |S|\right)$ and $n$. Hence $C_3$ depends on $\kappa, \gamma_0, R, \alpha,\beta$ and $\sup_{R_0 \leq b\leq R}\left(\frac{4}{r^2}\sup_{b<r} |S|\right)$.  Therefore the conclusion of the theorem follows.

\end{proof}




\section{The counting function on star-shaped Euclidean domains }\label{section: Euclidean domains}
 Let $\Omega \subset \R^n$ be a bounded domain in $\R^n$ containing the origin as an interior point such that there exists a a positive smooth function $\phi$ on the unit sphere such that 
\begin{align}\label{def of domain}
\Omega&=\{ x\in R^n: |x|<\phi(\frac{x}{|x|})\}.
\end{align}
 Let $(\rho, y)$ with $y\in \mathbb{S}^{n-1}$ be the polar coordinates of $\R^n$. Then $\partial \Omega = \{\rho = \phi(y)\}$. Note $\Omega$ is a star-shaped with respect to the origin. We will apply Theorem \ref{thm: estimate of spectral counting function - general case} to the counting function $N(\sigma)$ on $\Omega$.
 
We choose the potential function $f$ as
\begin{equation}
f(x) = \begin{cases} 
\frac{R^2 \rho^2}{4\phi^2(y)}, \quad, &x=(\rho, y)\neq 0;\\
0\quad, &x=0
\end{cases}
\end{equation} for some constant $R$. Then $b=2\sqrt{f}=\frac{R\rho}{\phi}$ at $x\neq 0$,  $\partial \Omega=\{b=R\}$ and $\Omega=\{b<R\}$.

We notice that the multiplicity of Steklov eigenvalue are invariant under scaling of the domain. Choose 
\begin{align}\label{choice of R}
R&=(\inf_{\mathbb{S}^{n-1}}\phi)\left(1+\sup_{\mathbb{S}^{n-1}}\phi^{-2}|\n^{\mathbb{S}^{n-1}} \phi|\right)^{-1/2}.
\end{align} Then by scaling of the domain, we can always have a large $R>8$. We choose $R_0 = 1$. Then Assumption \ref{basic assumptions} holds.

Next we will show that $\sup_{R_0 \leq r \leq R}\left(\frac{4}{r^2}\sup_{b<r}|S|\right)<1$.  Let $Y_1, Y_2, ..., Y_{n-1}$ be an orthonormal frame on the unit sphere $\mathbb{S}^{n-1}$ and denote the first and second order covariant derivatives of $\phi$ in this frame as $\phi_i, \phi_{ij}$, $i,j = 1, 2,..., n-1$. Choose $e_i = \rho^{-1} Y_i$ for $i = 1, 2,..., n-1$, and $e_n = \pd_\rho$ such that $\{ e_1, e_2, ..., e_n\}$ forms an orthonormal frame on $\R^n$. Then 
\begin{align*}
&(\nabla_{e_i} e_j)^\perp = - \rho^{-1} \delta_{ij} e_n, \quad \n_{e_i}e_n = \rho^{-1} e_i, \quad 1 \leq i , j \leq n-1.
\end{align*}
Let $\n_{e_i}f_i, \n_{e_i}\n_{e_j}f$ denote the first and second order covariant derivatives of $f$ in the frame $e_1, ..., e_n$. Then at $x\neq 0$, it follows that
\[
\n_{e_i} (f) = \begin{cases} 
-\frac{1}{2} R^2 \rho \phi^{-3} \phi_i, & i = 1,2,..., n-1; \\
\frac{1}{2} R^2 \rho \phi^{-2}, & i = n. 
\end{cases}
\]
\[
\n_{e_i} \n_{e_j} f = e_i (e_j (f) ) - \n_{e_i} e_j (f) = \begin{cases}
\frac{R^2}{2} (-\phi^{-3} \phi_{ij} + 3 \phi^{-4} \phi_i \phi_j) + \frac{R^2}{2} \phi^{-2}\delta_{ij}, & 1 \leq i, j \leq n-1; \\
 - R^2 \phi^{-3} \phi_i, & j = n, 1\leq i \leq n-1;\\
 \frac{R^2}{2 \phi^2}, & i = j = n.
\end{cases}
\]
Hence
\[
S(x) = f - |\n f|^2 = f (1 - R^2 \phi^{-4} |\n^{\mathbb{S}^{n-1}} \phi|^2 - R^2 \phi^{-2}), \quad x\neq 0.
\] Together with \eqref{choice of R}, this implies that $0<S<f$. Then we have $\sup_{r\leq R}\left(\frac{4}{r^2}\sup_{b<r}|S|\right)<1$.

Calculations above also implies that 
\[
(\sup_{\mathbb{S}^{n-1}} \phi)^{-2} (\inf_{\mathbb{S}^{n-1}}\phi)^2\left(1+\sup_{\mathbb{S}^{n-1}}\phi^{-2}|\n^{\mathbb{S}^{n-1}} \phi|\right)^{-1}\leq |\n b| \leq 1.
\] Here the left hand side is invariant under scaling.
And for all $1 \leq i , j \leq n$, we have
\begin{align*}
|\n_{e_i} \n_{e_j} f| \leq \frac{1+ 4\left(\sup_{\mathbb{S}^{n-1}}\phi^{-1}|\n^{\mathbb{S}^{n-1}}\n^{\mathbb{S}^{n-1}} \phi|+\sup_{\mathbb{S}^{n-1}}\phi^{-2}|\n^{\mathbb{S}^{n-1}} \phi|^2\right)}{2 \left(1+\sup_{\mathbb{S}^{n-1}}\phi^{-2}|\n^{\mathbb{S}^{n-1}} \phi|^2\right)}
\end{align*} with the right hand side invariant under scaling. Hence the tensor $T = \frac{1}{2} \delta - \n\n f$ is also bounded by a constant depending only on $\phi$ and invariant under scaling.

On $\R^n$ we can take the mean value constant $C_M = 1$. We also notice that the level sets of $b$ are homothetic which implies that
\[
Area(\{ b = s\}) = \frac{s}{R} Area(\partial \Omega), \quad r < R.
\]
Thus applying Theorem \ref{thm: estimate of spectral counting function - general case} implies the following result.
\begin{cor} For a bounded  domain $\Omega \subset \R^n$ satisfying \eqref{def of domain} for a smooth positive function $\phi$, the counting function of Steklov eigenvalue satisfies 
\[
{N}(\sigma) \leq C(\phi) Area(\pd \Omega)\left( \frac{1}{\inf_{\mathbb{S}^{n-1}} \phi} + \sigma \right)^{n-1}.
\]
\end{cor}

\section{The counting function and eigenfunctions on manifolds with boundary}
\label{section: manifolds with smooth boundary}

Consider a smooth compact manifold $(M^n, g)$ with smooth boundary $\pd M$. Instead of applying Theorem \ref{thm: estimate of spectral counting function - general case} directly, we will refine the estimates  in Section \ref{section: frequency} and Section \ref{section: counting function in general case} to obtain a more explicit upper bound of ${N}(\sigma)$ in Theorem \ref{intro - thm: main theorem}.  We will also study the decay behavior of eigenfunctions near the boundary and prove the Theorem \ref{intro - thm: behavior of eigenfunction near boundary}. 

Define the normal injective radius of the boundary (denoted by $\rm{nir}(\partial M^n)$) as the maximum of the set of radius $r>0$ such that the exponential map $\exp(x,t)=\exp_{x}(t \nu): \partial M \times [0, r)\rightarrow M$ is a diffeomorphism onto its image. Here $\nu$ is the inward unit normal of $\partial M^n$. By the compactness, this $\rm{nir}(\partial M^n)$ exists and is positive.

 Let $\rho(x) = d(x, \pd M^n)$ be the distance function to the boundary. Then $\rho$ is smooth in the $\rm{nir}(\partial M^n)$-neighborhood of $\partial M$ with $\rho = 0, |\n \rho| = 1$ on $\partial M^n$. And the metric $g$ can be written as 
\[
g = d\rho^2 + g_\rho,  \quad\rho \in [0, \rm{nir}(\partial M^n))
\]
where $g_\rho$ denote the induced metric on the level set of $\rho$ which is denoted as $\Sigma_\rho$. Then $\Sigma_0=\partial M$.

\begin{remark}
We assume the smoothness on $(M^n, g)$ and $\partial M^n$ for the convenience. Our argument works as long as $\rho(x) = d(x, \pd M)$ is $C^2$.
\end{remark}

\subsection{Upper bound of ${N}(\sigma)$: proof of Theorem \ref{intro - thm: main theorem}.}
 
\begin{proof}
Let $i_0>0$ be a constant slightly smaller than the normal injective radius. Choose $R =  k_0 i_0$ for some $k_0 > 1$ such that $R_0=R-i_0\geq 1$. We define the potential function \[f = \frac{1}{4}(R- \rho)^2\] for $0 \leq \rho \leq i_0$ and extend $f$ smoothly to the domain where $\rho (x)> i_0$ such that $f$ is well-defined on $M^n$. Then \[b=2\sqrt{f}= R- \rho\] for $0 \leq \rho \leq i_0$. And Assumption \ref{basic assumptions} holds. In particular, level sets of $b$ are just level sets of the distance function $\rho$ by $\{b=r\}=\Sigma_{\rho}$ for $\rho=R-r$. Through this section, we will use $\Sigma_{\rho}$ instead of $\{b=r\}$ in the argument. 
 
For $0 < \rho < i_0$, by $\n f = - \frac{1}{2} (R- \rho) \n \rho$, we have
\begin{align}
|\n b| &= 1, \quad S = f - |\n f|^2 = 0.
\end{align}
To calculate the tensor $T$, choose $\{e_1, e_2, ..., e_{n-1}\}$ as an orthonormal frame on level sets $\Sigma_\rho$ and $e_n = \nabla \rho$. With respect to $\{e_i\}$, there is $\rho_{ij} = h_{ij}(\Sigma_\rho), 1 \leq i, j \leq n-1$ where $\mathrm{I\!I}(e_i,e_j)=h_{ij}e_n$ as the second fundamental form on the level set $\Sigma_\rho$, and $\rho_{nn} = 0$. Denote the mean curvature on $\Sigma_\rho$ by $H$. Then for $\{0 < \rho < i_0\}$, by
 $\n\n f = - \frac{1}{2} (R - \rho) \n\n \rho + \frac{1}{2}\n \rho \otimes \n \rho$, we have 
\begin{align}\label{eqn: T}
T_{ij} &= \frac{1}{2} \delta_{ij} - \frac{1}{2} (R- \rho) h_{ij}(\Sigma_\rho),  \quad T_{ni} =0, \quad  T_{nn} = 0, \quad 1 \leq i,j\leq n-1.
\end{align} 

At first, we will refine key estimates of harmonic functions in Section 2. Let $u$ be any harmonic function well-defined on $\{R_0\leq b\leq R\}=\{0\leq \rho\leq i_0\}$. We start with the lower bound of $D^{\prime}(r)$ in Lemma \ref{lem: iterating the first variation formula}. Since $S=0$, the second term in the formula of $D^{\prime}(r)$ \eqref{deri of D2} vanishes. For the third term in \eqref{deri of D2}, by \eqref{eqn: T}, we have
\begin{align*}
\begin{split}
trT |\n u|^2 - 2 T(\n u, \n u) = & \left( \frac{n-1}{2} - \frac{R -\rho}{2} H\right) |\n u|^2 - \sum_{i = 1}^{n-1} u_i^2 +(R - \rho) \sum_{i,j = 1}^{n-1} h_{ij} u_i u_j  \\
& \leq \left( \frac{n-1}{2} + (R-\rho)(\frac{1}{2}C_H  + C_{\mathrm{I\!I}})\right) |\n u|^2.
\end{split}
\end{align*}with
\begin{align*}
C_H=\sup_{\rho \leq \rm{nir}(\partial M)}\sup_{\Sigma_{\rho}} |H|,\quad C_{\mathrm{I\!I}}=\sup_{\rho \leq \rm{nir}(\partial M)}\sup_{\Sigma_{\rho}}(|\mathrm{I\!I}|).
\end{align*}
Plug this inequality into \eqref{deri of D3}. Then we get
\begin{align}\label{eqn: differential inequality of D - boundary defining function}
\begin{split}
D'(r) 
&\geq 2r^{-1}D^2(r)I^{-1}(r) +\left(-\frac{n-1}{r}-C_H -2 C_{\mathrm{I\!I}}\right) D(r).
\end{split}
\end{align}

At the same time, plugging $S=0$ and \eqref{eqn: T} into the formula of $I^{\prime}(r)$ and $(\ln I(r))^{\prime}$ in Lemma \ref{lem: derivative of I}, we get refined bounds of $I^{\prime}(r)$ and $(\ln I(r))^{\prime}$ :
\begin{align}
I'(r) 
&\leq \frac{2D(r)}{r}-\left( \frac{n-1}{r}  - C_H\right) I(r),\label{eqn: differential inequality of I - boundary defining function}\\
- \frac{n-1}{r} -C_H &\leq ( \ln I(r) )' \leq \frac{2U(r)}{r}- \frac{n-1}{r} + C_H.\label{eqn: differential inequality of I' - boundary defining function}
\end{align}

Now apply both of \eqref{eqn: differential inequality of D - boundary defining function} and \eqref{eqn: differential inequality of I - boundary defining function} to the proof of Lemma \ref{prop: differential inequality for U}. We get a refined lower bound of $U^{\prime}(r)$:
\begin{align}\label{eqn: new lower bd of U'}
\begin{split}
U'(r) 
&\geq - 2(C_H+C_{\mathrm{I\!I}} )U(r)
\end{split}
\end{align}
Apply \eqref{eqn: new lower bd of U'} and \eqref{eqn: differential inequality of I' - boundary defining function} into the proof of Theorem \ref{thm: upper bd UI}. We obtain the refined growth estimates of $U(r)$ and $I(r)$ as:
\begin{align*}
U(r) &\leq U(R) e^{ 2(C_H+C_{\mathrm{I\!I}} )(R-r)}.
\end{align*} and for $R_0\leq r_2<r_1\leq R$
\begin{align}\label{eqn: new growth I}
e^{-C_H (r_1 - r_2)} \left(\frac{r_1}{r_2}\right)^{-n+1} \leq \frac{I(r_1)}{ I(r_2)} \leq (\frac{r_1}{r_2}) ^{2 U(R) e^{2(C_H+C_{\mathrm{I\!I}} )(R - r_2)} + (-n+1)} e^{C_H (r_1 - r_2)}.
\end{align} 

Next we run the same argument as the proof of Lemma \ref{L2 estimate on ball} using the refined growth estimate \eqref{eqn: new growth I}. Then the key inequality \eqref{eqn: doubling estimate on annulus regions} becomes
\begin{align*}
J(\lambda^{-3} R, R) &\leq ( \lambda^{- 1} e^{C_H R (\lambda - 1)} + 1 + \lambda ^{2 \sigma R e^{2(C_H+C_{\mathrm{I\!I}} )R(1 - \lambda^{-3})} +1} e^{C_H R (\lambda - 1)} ) J( \lambda^{-2} R, \lambda^{-1} R ).
\end{align*} and we obtain the refined $L^2$ estimate of $u$ on small balls as
\begin{align}\label{eqn: refinedL2}
 \int_{B(x, \delta)} u^2 &\leq \left(\lambda^{-1} e^{C_H(\lambda-1)R\lambda^{-3}}+ 1 + \lambda^{2 \sigma R e^{2(C_H+C_{\mathrm{I\!I}} )R(1 - \lambda^{-3})}+1}e^{C_H(\lambda-1)R\lambda^{-2}}\right)\int_{\lambda^{-2}R < b < \lambda^{-1}R} u^2
\end{align}  
with $\delta=\frac{1}{2}R (1-\lambda^{-1})\lambda^{-2}$ and $x\in [\lambda^{-2}R, \lambda^{-1}R]$. 

Now we can run the same argument as in the proof of Theorem \ref{thm: estimate of spectral counting function - general case} by
applying  \eqref{eqn: refinedL2} to \eqref{eqn: F2}. We pick a large $k_0 = \max(\frac{8}{7}, 1+i_0^{-1})$. Then
\begin{align*}
 R&=k_0i_0>i_0,\\
\alpha&=\frac{R_0}{R}=1-\frac{1}{k_0}>\frac{1}{8},\\
\lambda&=1+\frac{1}{(\alpha^{-1/3}-1)^{-1}+R\sigma}<\alpha^{-1/3}<2,\\
\delta&=\frac{R(1-\lambda^{-1})}{2}\lambda^{-2}=\frac{(\lambda-1)R}{2\lambda^{3}} < \frac{1}{2} \min \{i_0, \sigma^{-1}\},
\end{align*}
and similarly $R(1-\lambda^{-3} ) < 3 \min \{i_0, \sigma^{-1}\}$.
It follows that
 \begin{align}\label{eqn:new N}
 N(\sigma)
&\leq C(n)C_MAe^{2C_H \min\{i_0, \sigma^{-1} \}+ 2e^{6(C_H+C_{\mathrm{I\!I}} )\min\{i_0, \sigma^{-1} \}}} v_0^{-1}\left(\frac{1 }{i_0}+\sigma\right)^{n-1} \\
& \leq C(n)C_MAe^{ 4 \exp( 6(C_H+C_{\mathrm{I\!I}} )\min\{i_0, \sigma^{-1} \})} v_0^{-1}\left(\frac{1 }{i_0}+\sigma\right)^{n-1}.
\end{align} 
Here $C_M=C(n)e^{ C(n) i_0 \sqrt{K}}$ is the mean value constant and $K$ is from the Ricci lower bound $Ric \geq -(n-1)K$ in the $\rm{nir}(\partial M^n)$-neighborhood of $\partial M^n$. And $A=\sup_{0<\rho\leq i_0}A(\{\Sigma_{\rho}\})$ is the supremum of the area of level sets of $\rho$. 
We notice that by the first variation formula for the area, there is
\begin{align*}
A= \sup_{0<\rho\leq i_0}A(\{\Sigma_{\rho}\}) \leq {\rm Vol}(\partial M^n) e^{C_H i_0}.
\end{align*}
Plug this inequality into \eqref{eqn:new N}. Then it follows that \begin{align*}
N(\sigma)&\leq C(n)e^{C(n)i_0 \sqrt{K} +C_Hi_0 + 4e^{6(C_H+C_{\mathrm{I\!I}} ) \min\{i_0, \sigma^{-1}\} }}{\rm Vol}(\partial M^n) v_0^{-1} \left(\frac{1 }{i_0}+\sigma\right)^{n-1}.\\
\end{align*} 
Let $i_0 \to \rm{nir}(\partial M^n)$ and the conclusion \eqref{eqn:main thm1} and \eqref{eqn:main thm2} will follow.

\end{proof}
One immediate consequence of the upper bound of $N(\sigma)$ is the  following lower bound of Steklov eigenvalues.

\begin{cor}\label{cor: lower bound for eigenvalues}
Consider a smooth compact $(M^n, g)$ with smooth boundary $\partial M^n$. Let $0< \sigma_1 \leq \sigma_2 \leq ... \leq \sigma_j\cdots\nearrow +\infty$ be all positive Steklov eigenvalues on $(M^n, g)$. Then we have for all $j \geq 1$
\begin{align*}
\sigma_j&\geq C(n)e^{-\left(C(n) \sqrt{K}+C_H\right)\rm{nir}(\partial M)-4e^{6(C_H+C_{\mathrm{I\!I}} ) \min\{ {\rm nir}(\pd M), \sigma_1^{-1} \} }}v_0^{\frac{1}{n-1}}\left(\frac{j}{\rm{Vol}(\partial M)}\right)^{\frac{1}{n-1}}-\frac{1}{\rm{nir}(\partial M)}.
\end{align*} In particular, when $\sigma_j \geq \frac{c}{\rm{nir}(\partial M)}$ for some positive constant $c$, there is
\begin{align*}
\sigma_j&\geq C(n)e^{-\left(C(n) \sqrt{K}+C_H\right)\rm{nir}(\partial M) -4e^{6(C_H+C_{\mathrm{I\!I}} ) \min\{ {\rm nir}(\pd M), \sigma_1^{-1} \} }} v_0^{\frac{1}{n-1}}\left(\frac{j}{\rm{Vol}(\partial M)}\right)^{\frac{1}{n-1}}.
\end{align*}
\end{cor}
\begin{proof} Choose $\sigma=\sigma_j$ in Theorem \ref{intro - thm: main theorem}. Then we have
 \begin{align*}
j&=  N(\sigma_j)\\
&\leq C(n)e^{\left(C(n) \sqrt{K}+C_H\right)\rm{nir}(\partial M)+4e^{6(C_H+C_{\mathrm{I\!I}} ) \min\{ {\rm nir}(\pd M), \sigma_j^{-1} \}} }{\rm Vol}(\partial M^n) v_0^{-1}\left(\frac{1 }{{\rm nir} (\partial M^n)}+\sigma_j\right)^{n-1}.
\end{align*}
Taking the $(n-1)$th root of the above inequality and applying the fact $\sigma_j\geq \sigma_1$ yield the conclusions.

\end{proof}

\subsection{Behaviour of eigenfunctions near the boundary: proof of Theorem \ref{intro - thm: behavior of eigenfunction near boundary}}

\begin{proof}
Let $u$ be an Steklov eigenfunction corresponding to the Steklov eigenvalue $\sigma>0$ on $(M^n,g)$. Under the same setup of $f$ and $R$ as in the proof of \ref{intro - thm: main theorem}, we revisit Lemma \ref{lem: derivative of I} for $u$. It follows from $b=R-\rho$ and $|\nabla b|=1$ that
\begin{align*}
( \ln I(r) )' 
&= \frac{2U(r)}{r}-\frac{(n-1)}{r}+r^{1-n}I^{-1}(r)\int_{b=r}Hu^2ds.
\end{align*} This implies that
\begin{align*}
\frac{2U(r)-(n-1)}{r} + \inf H(\Sigma_\rho)  \leq (\ln I(r))' \leq \frac{2U(r)-(n-1)}{r}  + \sup H(\Sigma_\rho).
\end{align*}Integrating the above inequality on the interval $[r, R]$ yields 
\begin{align*}
&\left( \frac{R}{r}\right)^{- n+ 1}e^{(\int_r^R\inf H(\Sigma_{\rho})ds + \int_r^R\frac{u(s)}{s}ds)} \leq \frac{I(R)}{I(r)} \leq \left( \frac{R}{r}\right)^{ - n+ 1}e^{(\int_r^R\sup H(\Sigma_{\rho})ds + \int_r^R\frac{u(s)}{s}ds)}.\end{align*}

By the definition of $U(r)$ and the definition of Steklov eigenfunction in \eqref{eqn:Steklov1}, there is $U(R)=R\sigma$. Then, by the smoothness, we have in a sufficient small neighborhood of $\partial M=\{b=R\}$, i.e. $\rho \rightarrow 0$(or $r\rightarrow R$) that
\begin{align*}
U(r)&=U(R)+o(1)=R\sigma +o(1)\\
H(\Sigma_{\rho})&=H(\partial M)+o(1).
\end{align*} Plugging these approximation into the above inequality yields that\begin{align*}
\left( \frac{R}{r}\right)^{2 R \sigma + o(1) - n+ 1}e^{(R-r)\left(\inf H(\partial M) + o(1)\right)} \leq \frac{I(R)}{I(r)} \leq \left( \frac{R}{r}\right)^{2 R \sigma + o(1) - n+ 1}e^{(R-r)\left(\sup H(\partial M) + o(1)\right)},
\end{align*}
which can be equivalently written as 
\begin{align*}
\left( 1 - \frac{\rho}{R} \right)^{2 R \sigma + o(1) - n+ 1} e^{-(\inf H(\Sigma) + o(1))\rho} \leq \frac{I(r)}{I(R) } \leq \left( 1 - \frac{\rho}{R} \right)^{2 R \sigma + o(1) - n+ 1} e^{-(\sup H(\Sigma) + o(1))\rho}.
\end{align*}
Notice that $\left( 1 - \frac{\rho}{R} \right)^{R}\rightarrow e^{-\rho}$ as $R \to \infty$. Then the conclusion follows from taking $R \to \infty$ on above inequality.
\end{proof}

\section{ $C^{1, \alpha}$ domains in Riemannian manifolds}\label{section: c1alpha domains}

Let $\Omega$ be a bounded $C^{1, \alpha}$ domain in a smooth compact Riemannian manifold $(M^n, g)$. That is, every point in the boundary $\pd \Omega$ has a neighborhood where the boundary is the graph of some $C^{1,\alpha}$ function. In this case, the distance function to the boundary $\rho(x) = d(x, \pd \Omega)$ is only $C^{1,\alpha}$. Thus the method in Section \ref{section: manifolds with smooth boundary} does not apply directly. We need a regularized distance-like function. 

\begin{lem}\label{lem: regularized distance to the boundary}
 Let $\Omega$ be a bounded $C^{1, \alpha}$ domain in a smooth and comlete Riemannian manifold. For any $1>\epsilon > 0$, there exists a function $\eta \in C^{1, \alpha}(\bar{\Omega})$, satisfying 
\begin{align}\label{c1a}
&(1-\epsilon) \rho < \eta< (1+\epsilon) \rho, \quad  1-\epsilon \leq |\n \eta| < 1+\epsilon, \quad and \quad  |\n\n \eta| \leq C \rho^{\alpha -1} 
\end{align}
in a neighborhood of the boundary, where the constants $C$ depends on $\epsilon$, $\alpha$, and the geometry of a tubular neighborhood of $\pd \Omega$ in $M.$
\end{lem}

\begin{proof} 
 
By the result of Anderson (see Main Lemma 2.2 in \cite{And90}), for any $1>\epsilon>0$, there exists a radius $\delta>0$ such that one has $C^{1, \alpha}$ harmonic coordinates on $B_{ \delta}(y)$ for any $y \in \pd \Omega$. This radius $\delta$ depends on $\epsilon, \alpha$, the Ricci curvature, and the injective radius. With respect to this harmonic coordinates on $B(y, \delta)$ denoted by $(x^1, x^2, ..., x^n)$, the Riemannian metric $g$ is $C^{1, \alpha}$-close to the Euclidean metric, i.e. $|g_{ij} - \delta_{ij}|_{C^{1, \alpha}} < \epsilon$.

 By the compactness, we can cover $\pd \Omega$ with finite small balls $B_{\delta}(y)$. We will construct a distance-like function satisfying \eqref{c1a} in each $B_{\delta}(y)$ and glue them up by a smooth partition of unity. Then we extend this function smoothly into $\Omega$ to get $\eta$. 

Fix one $B_{\delta}(y)$. Let $\phi$ be a nonnegative smooth function in a neighborhood of $B_{ \delta}(y)$ which has a compact support in $B_{ \delta}(y)$ and equals to $1$ in $B_{\delta/2}(y)$ and $0$ outside of $B_{ \delta}(y)$. Define a parabolic operator
\begin{align*}
L =\pd_t - \phi \Delta^E - (1-\phi) \Delta_g
\end{align*}
where $\Delta_g$ is the Laplacian operator with respect to the Riemannian metric $g$, and $\Delta^E$ is the Laplacian operator with respect to the Euclidean metric in these harmonic coordinates $(x^1, x^2, ..., x^n)$. Note that $L$ is the standard Euclidean heat operator on $B_{\delta/2}(y)$. Then there exits a $C^{1,\alpha}(\bar{\Omega})$ solution of the initial-boundary value problem 
\[
\begin{cases}
L u(x, t) = 0, & \text{in }B_{ \delta}(y)\cap \Omega; \\
u(x, t) = 0, & \text{on } B_{ \delta}(y)\cap \pd \Omega; \\
u(x, 0) = \phi(x)\rho(x). 
\end{cases}
\]
We can take $t_0> 0$ small enough such that 
\begin{align*}
&(1- \epsilon) \rho(x) < u(x, t) < (1+\epsilon ) \rho(x), \quad and \quad 1-\epsilon < |\n u(x, t)| \leq 1 + \epsilon,
\end{align*}
for all $t \in [0, t_0]$ and $x \in B_{\delta/2} (y) \cap \Omega$.

Next we derive the estimate of $\nabla \nabla u$ near $y$. For any constant unit vector $\vec{v}$ and sufficiently small $h > 0$, define $D^h \pd_i u(x,t) = \pd_i u(x+h \vec{v}, t) - \pd_i u(x, t)$. Since $L$ has constant coefficients in $B_{\delta/2}(y)$, the function $D^h \pd_i u(x,t)$ is also a solution of the heat equation:
\begin{align*}
L D^h \pd_i u(x,t)& = 0, \quad x \in B_{\delta/2}(y)\cap \Omega.
\end{align*}
Note that $|D^h \pd_i u(x,t) | < C h^\alpha$, by parabolic gradient estimate we have 
\begin{align*}
|\n D^h \pd_i u|(x, t)& \leq C h^{\alpha} \rho(x)^{-1} \wedge t^{-1/2}, \quad (x, t) \in B_{\delta / 4}(y) \cap \Omega \times (0, t_0].
\end{align*}
For any fixed $t\in (0, t_0]$, pick $\rho(x) < t^{1/2}$ and $h = \rho(x)/2$. Then it follows that 
\begin{align*}
|\n\n u|(x + h \vec{v}, t)& \leq | \n\n u|(x, t) + C \rho(x)^{\alpha - 1}.
\end{align*}
By iterating the above inequality, we can find a sequence of points $x_0 = x, x_1, x_2, ..., $ in $B_{\delta/2}(y) \cap \Omega$ along the geodesic joining $x$ to the boundary $\pd \Omega$ such that $\rho(x_i) = 2^{-i} \rho(x)$ and 
\[
\begin{split}
|\n\n u|(x_i, t) \leq & |\n\n u|(x_0, t) + C \sum_{j = 0}^i (2^{-j} \rho(x))^{\alpha - 1}\\
\leq & |\n\n u|(x_0, t) + C 2^{(i+1)(1-\alpha)}\rho(x)^{\alpha-1} \leq C \rho(x_i)^{\alpha-1}.
\end{split}
\]
Here the constant $C$ varies from line to line. It is not hard to see that $x_i$ can be made arbitrary, hence we have established the desired second order estimate for $u$ is a smaller neighborhood around $y$. Then we choose $u(x, t_y)$ for some $t_y<t_0$ as the desired distance-like function in $B_{\delta}(y)$.

\end{proof}

Now we are ready to prove Theorem \ref{intro-thm: estimate of counting function - C1alpha domain} in the introduction. 

\begin{proof}[Proof of Theorem \ref{intro-thm: estimate of counting function - C1alpha domain} ]
Fix a positive constant $\epsilon \in (0,1)$ and $R > 1$. Choose $R_0= R - i_0$ for  a sufficiently small $i_0> 0$. Define $f = \frac{1}{4}(R-(1-\epsilon)\eta)^2$ and $b = 2 \sqrt{f} = R - (1-\epsilon)\eta$, where $\eta$ is the function from Lemma \ref{lem: regularized distance to the boundary}. Then 
\begin{align*}
&\n f = \frac{1}{2} (R - (1-\epsilon) \eta) (- (1-\epsilon) \n \eta), \quad |\n f|^2 = (1-\epsilon)^2 |\n \eta|^2 f, \\
&\n\n f = \frac{(1-\epsilon)^2}{2} \n \eta \otimes \n \eta  - \frac{1-\epsilon}{2} b \n\n \eta.
\end{align*}
By Lemma \ref{lem: regularized distance to the boundary}, we can check that Assumption \ref{basic assumptions} holds. Moreover we have
\begin{align*}
&0 \leq S = f -|\n f|^2 \leq (1 - (1-\epsilon)^2) f, \quad |\n\n b|\leq C \rho^{\alpha - 1}, \quad |T|\leq C \rho^{\alpha - 1} 
\end{align*}
in the domain $\{R_0 < b < R\}$. By taking $i_0$ sufficiently small, the condition $4\|S\|_{L^\infty(b< r)} < r^2$ is satisfied. Thus we can apply Theorem \ref{thm: estimate of spectral counting function - general case} to get the conclusion.
\end{proof}
We also obtain the following result about the behavior of Steklov eigenfunctions near a $C^{1, \alpha}$ boundary. 
\begin{thm}\label{thm:C1alpha}
Let $\Omega$ be a bounded $C^{1, \alpha}$ domain in a smooth compact Riemannian manifold $(M^n, g)$. Let $\rho$ be the distance-to-boundary  function and denote its level sets as $\Sigma_\rho$. For any Steklov eigenfunction $u$ with respect to the eigenvalue $\sigma>0$, there is
\begin{align*}
 &e^{ -C\rho^\alpha- 2 (1+\epsilon) \sigma \rho} \leq \frac{\int_{\Sigma_\rho}  u^2 }{ \int_\Sigma u^2 } \leq   e^{ C\rho^\alpha - 2(1-\epsilon) \sigma \rho}
\end{align*}
when $\rho$ is sufficiently small. Here $C$ is a constant depending on $\alpha$ and the geometry near the boundary.
\end{thm}
\begin{proof}
For any $1>\epsilon > 0$, let $\eta$ be the function from Lemma \ref{lem: regularized distance to the boundary}. Define the potential function $f$ and $b$ as in the proof of Theorem \ref{intro-thm: estimate of counting function - C1alpha domain}. Then we have 
\begin{align*}
&(1-\epsilon) R \sigma \leq U(R) \leq (1+\epsilon) R \sigma
\end{align*}
in a small neighborhood of $\pd \Omega$. 

By Lemma \ref{lem: derivative of I} and Lemma \ref{lem: regularized distance to the boundary}, we have
\begin{align*}
&\frac{2(1-\epsilon) R \sigma }{ r } - \frac{C}{  (R-r)^{1-\alpha}} \leq (\ln I(r))' \leq \frac{2 (1+ \epsilon) R \sigma }{ r } + \frac{C}{  (R-r)^{1-\alpha}} .
\end{align*}
Integrating the above inequality on the interval $[r, R]$ yields 
\begin{align*}
&\left( \frac{R}{r}\right)^{2 (1-\epsilon) R \sigma} \exp \left( -\int_r^R \frac{C}{(R-s)^{1-\alpha}} ds\right) \leq \frac{I(R)}{I(r)} \leq \left( \frac{R}{r}\right)^{2 (1+\epsilon) R \sigma} \exp \left( \int_r^R \frac{C}{(R-s)^{1-\alpha}} ds\right)
\end{align*}
which implies that
\begin{align*}
\begin{split}
\left( 1 - \frac{(1+\epsilon) \rho}{R} \right)^{2 (1+\epsilon) R \sigma}  \exp \left( - \int_r^R \frac{C}{(R-s)^{1-\alpha}} ds\right)  \leq \frac{I(r)}{I(R) } \\
 \leq  \left( 1 - \frac{(1-\epsilon) \rho}{R} \right)^{2 (1-\epsilon) R \sigma}  \exp \left( \int_r^R \frac{C}{(R-s)^{1-\alpha}} ds\right)
\end{split}
\end{align*}
Let $R \to \infty$ and the proof is finished by adjusting $\epsilon$. 
\end{proof}

\end{document}